\documentclass[12pt]{elsarticle}

\usepackage{amssymb}
\usepackage{caption, color, url}

\begin{document}

\begin{frontmatter}

\title{A Matheuristic for the Electric Vehicle Routing Problem with Time Windows}

\author[1]{Maurizio Bruglieri\corref{cor1}}
\author[2]{Ferdinando Pezzella}
\author[3]{Ornella Pisacane}
\author[2]{Stefano Suraci}

\cortext[cor1]{maurizio.bruglieri@polimi.it}

\address[1]{Dipartimento di Design, Politecnico di Milano, Milano, Italy}
\address[2]{Dipartimento di Ingegneria dell'Informazione, Universit\'{a} Politecnica delle Marche, Ancona, Italy}
\address[3]{Facolt\'{a} di Ingegneria, Universit\'{a} degli Studi e-Campus,
Novedrate (Como), Italy}

\begin{abstract}
The main goal of this paper is to time-effectively route and schedule a fleet of Electric Vehicles (EVs) on a road network in order to serve a set of customers. In particular, we aim to propose an optimized route planning by exploiting the advantages of these vehicles. Nowadays, in fact, electromobility plays a key role for reducing the harmful emissions due, instead, to the use of traditional vehicles. The starting point of this research is represented by the fact that the advanced recent technologies for the EVs allow also partially recharging their batteries. In this work, an Electric Vehicle Routing Problem with Time Windows (E-VRPTW) is addressed from a time effective point of view under the assumption that partial recharges are also allowed. For this purpose, the E-VRPTW is mathematically formulated as a Mixed Integer Linear Program in which both the total number of EVs used and the total time spent by them outside the depot are minimized. Due to the NP-hardness of the problem, a Variable Neighborhood Search Branching (VNSB) matheuristic is also designed for determining good quality solutions in reasonable computational times. Numerical results carried out on some benchmark instances taken from the literature provide useful insights regarding both the solution quality of the proposed formulation, compared to a previous one, and the performances of the VNSB.

\end{abstract}

\begin{keyword}
Route Planning \sep Electric Vehicles \sep Partial Battery Recharge \sep Variable Neighborhood Search \sep Local Branching.

\end{keyword}

\end{frontmatter}


\section{Introduction}
\label{intro}
Nowadays, to provide eco-sustainable transportation solutions represents one of the most significant steps toward the design of smarter cities in which the attention is mainly devoted to citizens' quality of life.

To this end, the Electric Vehicles (EVs) play a key role. In fact, they are less noisy than the traditional Internal Combustion Engine Vehicles (ICEVs) and guarantee also no harmful emissions of $CO_2$. These advantages lead to an \textit{environmental sustainability}.

The European Commission has recently affirmed that the one-fifth of the EU's total emissions of $CO_2$ is due to road transport \cite{EU}. Moreover, it has also remarked that from $1990$ to $2007$, they decreased in other sectors while increased of about $36\%$ in the transport one. These considerations justify one of the main EU's targets by $2020$ such as their reduction of about $20\%$ with regard the levels reached in $1990$. For this purpose, several governments are actually adopting specific programs such as, for example, to limit the access of the ICEVs to the some urban areas. These areas are usually around the historical center in which the traffic of ICEVs can be forbidden either all day or in specific time periods. In this case, for example, the EVs can be used for providing particular categories of citizens (e.g., disables and elderlies) with efficient door-to-door transportation services also in these particular urban areas. This point remarks the \textit{social sustainability} of the EVs.

Although the purchase cost of an EV is higher than a conventional one, its use is by far economically advantageous. In fact, a recent comparison analysis \cite{Feng}, between a conventional diesel truck (Isuzu N-Series) and an electric engine truck (Navistar E-star), shows that the average operating cost of the former is $\$0.23/$miles against $\$0.09/$miles for the latter.
Moreover, the authors remark that electric trucks are $50\%$ less expensive to maintain than the conventional ones. This leads to an \textit{economical sustainability}.

These aspects also contribute to increase the interest of the automotive industry to produce more efficient EVs and to overcome some limitations of them. Among these, the poor battery range has to be considered. In fact, as noted in \cite{Feng}, a full battery recharge allows traveling at most $100$ miles, on average and this autonomy is usually reduced by $50\%$ due to the so called ``range anxiety''. This is a disadvantage in the cases in which the EVs are used for the long-distance transport, i.e., as commercial vehicles rather than  passenger ones.

In addition, the need of many stops during a trip usually increases the duration of the routes, aspect taken under control especially in the \textit{Distributive Logistics} that is receiving a strong impulse also thanks to the e-commerce (\textit{Amazon},\textit{ eBay}, to cite a few companies).

Thanks to the technological advancements, the modern electric batteries do not suffer from the so-called ``memory effect'' and then, it is possible to partially recharge the EVs during the trip with a consequent reduction of the times spent at the RSs.

The aim of this paper is to address a time-effective Electric Vehicle Routing Problem with Time Windows (E-VRPTW) which consists in finding a set of routes by EVs, starting/returning from/to a common depot, in order to handle a set of customers requiring a service inside specific time windows. Thus, the objective of this problem is the minimization of the total time spent, by the EVs used, outside the depot.

\begin{figure}[t!]
\begin{center}
\begin{tabular}{c}
\includegraphics[width=10 cm]{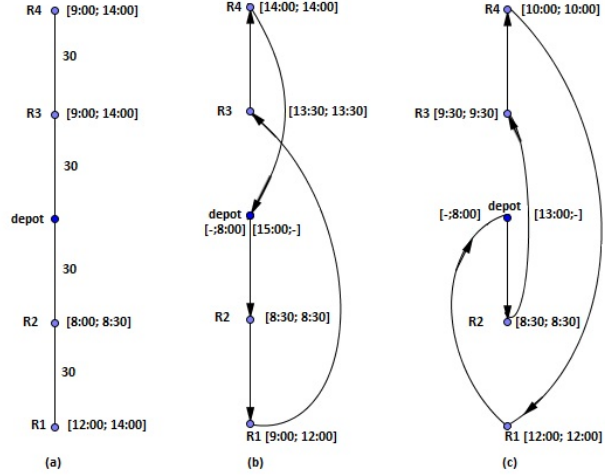} \\
\end{tabular}
\caption{A numerical example} \label{example}
\end{center}
\end{figure}

In order to remark the practical impact of such a new problem, a numerical example is shown in Figure \ref{example} (a) where it is assumed that four customer requests have to be handled. In particular, the requests R$1$ and R$4$ are distant $60$ Km from the depot while R$2$ and R$3$, $30$ Km. R$1$ has a time window equal to $[12$:$00;14$:$00]$, R$2$ equal to $[8$:00$;8$:$30]$ and finally, R$3$ and R$4$ equal to $[9$:$00;14$:$00]$. Finally, the vehicle speed is assumed unitary (i.e., $1$ Km per minute) and all the service times at the customers are supposed to be negligible.

In Figure \ref{example} (b), the solution that minimizes the total travel distance is shown while in Figure \ref{example} (c), the solution that minimizes the total route duration.
In both the solutions, the arrival and the departure times of the EV to and from each customer are indicated beside it, between square brackets, respectively. While, for the depot, two time intervals of this kind are used: the former related to the departure time and the latter, to the arrival time.

It is possible to observe that, in the first solution, the vehicle first arrives to the customer R$2$, with the strictest time window, at $8$:$30$. Then, it moves to the nearest customer R$1$ where it arrives at $9$:$00$. In this case, it has to wait three hours before starting the service since R$1$ cannot be served before $12$:$00$ and not later $14$:$00$. Then, it reaches R$3$ at $13$:$30$ and finally, R$4$ at $14$:$00$, returning to the depot at $15$:$00$. In this way, the total route duration is equal to seven hours with three hours of waiting time.

In the second solution, instead, the vehicle arrives to R$2$ at $8$:$30$ and then, it moves to R$3$ that is distant $60$ Km and that can be handled immediately. Then, it goes to R$4$ and finally, to R$1$. Thus, the EV can return to the depot at $13$:$00$. In this case, the total route duration is equal to five hours with zero waiting time.

Comparing the two solutions, it can be observed that the total travel distance of the first solution is equal to $240$ Km while, the one of the second solution is equal to $300$ Km. However, especially in the sector of the last mile Logistics, a saving of $60$ Km may be considered negligible with regard to a saving of $2$ hours, as guaranteed by the second solution.

Due to the poor battery range, an EV may require to be recharged more than once during its trip and therefore, intermediate stops at the RSs have been also scheduled. Obviously, this may cause an increment of the total time duration.
In order to take under control that increment, in this paper, an EV is allowed being also partially recharged at a RS.

Therefore, the final solution aims firstly to minimize the number of EVs used and then, the total time spent by the EVs outside the depot, i.e., for recharging, traveling and waiting. Such a time is obtained as the sum of the differences of the EVs' arrival times to the depot with their starting times from the depot.

For addressing such a new optimization problem, a Mixed Integer Linear Programming (MILP) formulation under both time and battery constraints is proposed.

Moreover, to tackle the problem, despite its NP-hardness, a matheuristic, i.e., a Variable Neighbor Search Branching (VNSB) \cite{Hansen}, is also designed.

A preliminary version of this work has been already presented in \cite{Bruglieri} where computational experiments carried out on a sub-set of benchmark instances gave very promising results. However, although this work shares the methodological approach with the previous one (i.e., VNSB), here we propose a different MILP formulation. In fact, in order to optimize the total time spent by the EVs outside the depot, it simply minimizes the difference between the arrival time at the depot and the starting time from the depot. To this end, the depot is properly cloned into two distinct sets in order to properly differentiate the routes. The resulting new MILP formulation is defined with fewer decision variables and constraints than the previous one.

The rest of the paper is organized as in the following: Section \ref{literature} reviews the state of the art, Section \ref{model} describes the MILP formulation proposed for the problem under consideration while Section \ref{vns} presents the designed VNSB. Finally, Section \ref{results} shows computational results obtained on several benchmark instances and Section \ref{conclusions} concludes the work remarking also possible future developments.

\section{Literature review}
\label{literature}
The problem addressed in this paper aims to handle a set of customers within specific time windows by using a fleet of EVs. Due to the general features of such a problem, it can be seen as a variant of the well-known Vehicle Routing Problem (VRP). The aim of a VRP is to find a set of routes, one for each vehicle of the fleet, in order to handle a set of requests, optimizing one or more performance measures (e.g., minimizing the total travel cost). Each route has to start from a depot (that could be common to all the vehicles) and has to return to a depot (that could be the same of the initial one). According to the number of depots, the VRPs are classified into the Single Depot VRPs (SDVRPs) and the Multi-Depot VRPs (MDVRPs). Moreover, in the case in which the vehicles have a limited capacity, the VRPs become capacitated (i.e., CVRP). It is worth noting that the vehicles could have either all the same capacity (i.e., homogeneous fleet) or different capacity (i.e., heterogeneous fleet). In addition, if each customer request has to be satisfied within a specific time window, a VRP with Time Windows (VRPTW) is addressed. The readers are referred to \cite{Kumar} and \cite{Laporte} for two recent surveys on the VRP and its variants.
\\Indeed, due to the recent developments of the automotive industry and to several environmental problems, new variants of the VRP have been proposed in the literature. Among these, the ones of our interest aim to control the harmful emissions of $CO_2$ and to manage fleets of alternative fuel vehicles. Such a problem is not new to literature. In \cite{Figliozzi}, the author aims to route a fleet of capacitated vehicles in order to control the harmful emissions: Emissions VRP with time windows. Therefore, a cost associated with the $CO_2$ emissions is introduced. Then, a first Mixed Integer Linear Programming (MILP) formulation is proposed for minimizing a total cost due in part to the vehicles used, to the travel distance, to the route duration and then, to the emissions. Finally, a second MILP formulation is derived in which priority is given to the minimization of the vehicles used, then to the emissions and finally, to both the total distance and the route duration. It is worth noting that both the MILP formulations impose, together with the traditional VRP constraints, also feasibility conditions for the vehicle capacity and the customer time windows. For such a problem, the author designs a heuristic approach where, in a first constructive phase, feasible routes are determined with the aim of minimizing the total number of vehicles used; in a improvement phase, the emissions are minimized.

In \cite{Erdogan}, a MILP formulation for a VRP with a fleet of fuel-powered vehicles is proposed: the so called Green VRP (G-VRP). The objective function, to minimize, denotes the total travel distance. Moreover, together with the traditional VRP constraints, a condition on the maximum duration of each route is also imposed. Once a vehicle arrives at a fueling station, its fuel level is set to the maximum capacity. On the contrary, the arrival at a customer implies a fuel consumption. It is worth noting that the authors do not model feasibility conditions on both the time windows associated with the customers and the vehicle capacity. From a methodological point of view, two heuristics are proposed: one modifies the well-known Clark and Wright savings algorithm while the other is a density-based clustering algorithm.

In \cite{Conrad}, instead, a first variant of the VRP, called the Recharging VRP, in which the vehicles of the fleet are allowed to be recharged during their trips is introduced. The proposed MILP formulation aims to firstly minimize the number of vehicles employed and then, a total cost expressed as the sum of the total travel cost, of the cost due to the route duration and to the vehicle recharges. The two main assumptions of this work concern both the recharging points and times. In particular, with reference to the former, the authors assume that each customer location is also a recharging station while, regarding the latter, only fixed recharging times are considered.

A literature contribution overcoming some of the limitations remarked above is the one proposed in \cite{Schneider}. The aim of this work is to address an E-VRPTW with homogeneous fleet and recharging stations. The problem is mathematically formulated as a MILP model, on a complete and direct graph, in which the set of vertices is represented by both the RSs and the customers. For modeling reasons, the RSs are cloned and, beyond the traditional VRPTW constraints, battery feasibility conditions are also imposed. The set of routes is found in order to handle all the customers in time, firstly minimizing the number of EVs employed for delivering the services and then, the total travel distance. The authors also generate a set of benchmark instances for the E-VRPTW from the one already described in \cite{Solomon} for the VRPTW. With the aim of providing good quality solutions in reasonable computational times, they also design a hybrid metaheuristic combining a variable neighborhood search with a tabu search. The computational results obtained by the authors through the hybrid meta-heuristic show the high performance of the proposed solution approach. However, an assumption of such a work consists in always fully recharging the EVs at the RSs.

On the other hand, the concept of the partial recharge has been already addressed in literature. For example, in \cite{Righini}, the authors address the G-VRP 
introducing different recharging technologies and also, the partial battery recharges. For this purpose, they propose a MILP formulation in which the objective function, to minimize, denotes the total cost due to both a fixed component (related to the use of the EVs) and a variable part (related to the battery recharges). Moreover, they design column generation algorithms for solving this problem. However, the authors do not take into consideration the time windows on the customer services.

In \cite{Alesiani}, the authors address the problem of efficiently routing a homogeneous fleet of EVs, with a single depot but without time windows constraints. Under the assumption that the EVs are always fully recharged, they aim to route the fleet in order to minimize both the total travel distance and the total waiting time of the drivers at the RSs. For this purpose, they limit the number of stops that a vehicle can do for recharging itself during the trip. Moreover, they also take into account the energy consumption/gain by including, for example, the losses due to both the rolling and the aerodynamic resistances. For this problem, the authors propose a MILP formulation together with a genetic algorithm for solving it on large-scale instances. In \cite{Wang}, the single-origin single destination VRP is firstly addressed with the aim of minimizing both the route time duration and the times spent for recharging the vehicle. For this problem, the authors propose a Mixed Integer NonLinear Programming formulation where the non-linearity is due to the need of modeling the energy consumption constraints.
By exploiting some general properties of such a problem, they address also the multi-vehicles variant. In \cite{Preis}, an E-VRPTW is addressed considering full battery recharges and the objective function, to minimize, denotes the total energy consumption. It is expressed in function of both the driving resistances and the loading weight. For this purpose, the authors formulate a MILP model and design a meta-heuristic tabu search. The computational results show promising percentages of energy saving during the routes. The same authors in \cite{Frank} propose other different MIP formulations in which the objective function, to minimize, considers the cost depending on the energy consumption of both the empty vehicles and the payloads. Finally, they design a column generation based solution approach for real-life alike instances. In \cite{Felipe}, several heuristic approaches are proposed with the aim of addressing a G-VRP in which the EVs can be recharged by using different technologies. Under the limitation that the RS at the depot can be used only during the night, constructive greedy heuristics, a deterministic local search and solution approaches based on the simulated annealing framework are proposed. According to the computational results, they conclude that the exhaustive local search performs better than the others on medium-scale instances while, the simulated annealing based heuristics give the best on large-scale instances. \\In \cite{Schneiderb}, a VRP with Intermediate Stops (VRPIS) is described for which the authors design an adaptive variable neighborhood search. It is also shown how the mathematical formulation of the VRPIS can be used for solving an E-VRP with recharging facilities but under the assumption of only full battery recharges. The objective function, to minimize, denotes a total cost due to both a variable part related to the routes and a fixed component related to the vehicles used. It is worth noting that time windows for the services at the customers are not considered but the arrival at each vertex is limited to be less or equal to the maximum duration allowed for each route.

In \cite{Montoya}, a G-VRP is addressed with a fleet of zero emission vehicles without planning intermediate stops, during the routes, for recharging reasons. The problem is addressed from a methodological point of view starting from the mathematical formulation proposed in \cite{Erdogan}. A two-steps heuristic is designed and tested where, in the first phase, a set of feasible routes is built through both route-first cluster-second heuristics and an ad hoc insertion procedure for the alternative fuel stations. In the second phase, instead, these feasible routes are assembled through a set-partitioning formulation in order to find solution to the G-VRP.

\section{A MILP formulation of the Time-Effective E-VRPTW}
\label{model}
Starting from the MILP formulation described in \cite{Schneider}, one of the contributions of our work is 
modeling the concept of partial recharges. This is mainly motivated by the following two reasons. On one hand, the immediate effect is a decrement of the total recharging time. On the other hand, partial recharges may also allow serving in time customers who could not be reached by the same driver within their time window, otherwise. In this way, the feasible region of our variant of the E-VRPTW is wider than the original one proposed by \cite{Schneider} and then, we could also obtain solutions with fewer routes or lower total traveled distance. 
\\The problem is represented on a directed graph $G=<V,A>$ where the set of the vertices $V=N\cup F\cup D_0 \cup D_{0'}$ contains the sets: $N$ of $n$ customers, $F$ of $f$ RSs plus $f'$ clones of them necessary to allow modeling several visits to the RSs as elementary routes, $D_0$ and $D_{0'}$ representing the sets of the depot clones used at the beginning of each route and at the end, respectively. The clones of the depot are necessary to properly model the ending time of each route.
In addition, the subset $V'$ denotes $V\setminus\{D_{0'}\}$ while $V''$ is $V\setminus\{D_0\}$.
 The set $A$ of arcs contains all the ordered pairs of vertices except those between depot clones as well as those from each RS to its clones. For each arc $(i,j) \in A$, both a travel distance and time, denoted by $d_{ij}$ and $t_{ij}$, respectively, are given. Since it is assumed that the fleet is homogeneous, the capacity of each EV is denoted by $C$ while the battery capacity by $Q$. The average EV speed is $v$ while the consumption and the recharging rate of the battery is $r$ and $g$, respectively.
\\The request of a customer $i \in N$ is expressed in function of a known demanded quantity $q_i \geq 0$ modeling a pickup while the related time window is indicated as $[e_i, l_i]$. Under the assumption of hard time windows constraints, the service cannot be performed before $e_i$ and later than $l_i$. It is worth noting that a time window $[e_0, l_0]$ is defined also for the depot, where $l_0$ indicates the maximum time allowed for coming back to depot. 
The service time of each customer $i$ is denoted as $s_i \geq 0$.

The E-VRPTW consists in finding a set of vehicle routes starting from a vertex of $D_0$ and ending to a vertex of $D_{0'}$ in such a way that each node of $N$ is visited from exactly one vehicle route within its time window, the capacity of the vehicle is not exceeded, and the battery level never becomes negative. To the aim of satisfying the latter condition, the routes can also pass to the nodes of $F$ obtaining a recharge, linearly proportional by $g$ to the recharging time, until at most $Q$.

The goal is to minimize firstly the number of vehicle routes and then, the total time spent by the vehicles outside the depot, i.e., the sum of the recharging, traveling and waiting times. Such a total time is obtained as the sum of the differences of the EVs' arrival times to the depot with their starting times from the depot.

This problem is mathematically formulated by introducing the decision variable, one for each arc $(i,j)\in A$, $x_{ij}$ equal to $1$ if $(i,j)$ is traversed, $0$ otherwise. In order to model the arrival time and remaining vehicle capacity at a vertex $i \in V$ as well as the battery level before leaving $i$, the following decision variables $\tau_i$, $u_i$ and $y_i$ are introduced, respectively. With the aim of considering also partial battery recharges, the decision variable $\xi_i$ denotes the battery level reached at the RS $i \in F$. In this way, the battery level reached at each RS is established during the optimization process.


The two-indices MILP formulation is detailed in the following:

\begin{equation}
\label{eq.objective_function}
\min \hspace{0.5mm} \sum_{i \in D_0, j \in N \cup F} x_{ij}+\hspace{0.5mm} \sum_{i \in D_0, j \in D_{0'}|i \neq j} \frac{(\tau_j-\tau_i)}{l_0|N||D_0|}
\end{equation}

\begin{equation}
\label{eq.1}
\sum_{j \in V''|j \neq i} x_{ij}=1 \hspace{3mm} \forall i \in N
\end{equation}

\begin{equation}
\label{eq.2}
\sum_{j \in V''|j \neq i} x_{ij}\leq 1 \hspace{3mm} \forall i \in F
\end{equation}

\begin{equation}
\label{eq.3DepIniziale}
\sum_{j \in N\cup F} x_{ij}\leq 1 \hspace{3mm} \forall i \in D_0
\end{equation}

\begin{equation}
\label{eq.3}
\sum_{i \in N\cup F} x_{ij}\leq 1 \hspace{3mm} \forall j \in D_{0'}
\end{equation}

\begin{equation}
\label{eq.4}
\sum_{i \in V''|j \neq i} x_{ji}- \sum_{i \in V'|j \neq i} x_{ij}=0 \hspace{3mm} \forall j \in N \cup F
\end{equation}

\begin{equation}
\label{eq.5_1}
\tau_i+(t_{ij}+s_i)x_{ij}-l_0(1-x_{ij})\leq \tau_j \hspace{3mm} \forall i \in N \cup \{D_0\}, \forall j \in N \cup \{D_{0'}\}| i \neq j
\end{equation}

\begin{equation}
\label{eq.6_1}
\tau_i+t_{ij}x_{ij}+g(\xi_i-y_i)-(l_0+gQ)(1-x_{ij})\leq \tau_j \hspace{3mm}  \forall i \in F, \forall j \in N\cup \{D_{0'}\}| i \neq j
\end{equation}

\begin{equation}
\label{eq.7}
e_j\leq \tau_j\leq l_j \hspace{3mm} \forall j \in V \setminus\{F\}
\end{equation}

\begin{equation}
\label{eq.9}
u_i \leq C\hspace{3mm}\forall i \in D_0
\end{equation}

\begin{equation}
\label{eq.10}
y_j \leq y_i-rd_{ij}x_{ij}+Q(1-x_{ij}) \hspace{3mm} \forall i \in N\cup \{D_0\}, \forall j \in V''| i \neq j
\end{equation}

\begin{equation}
\label{eq.11}
y_j\leq \xi_i-rd_{ij}x_{ij}+Q(1-x_{ij}) \hspace{3mm} \forall i \in F, \forall j \in V''| i \neq j
\end{equation}

\begin{equation}
\label{eq.12}
u_j \leq u_i - q_ix_{ij} + C(1 - x_{ij}) \hspace{3mm} \forall i \in V',\forall j \in V''|i \neq j
\end{equation}

\begin{equation}
\label{eq.13_bis}
\tau_j\leq l_0\sum_{i \in V'} x_{ij} \hspace{3mm} \forall j \in D_{0'}
\end{equation}

\begin{equation}
\label{eq.13_ter}
\tau_i\leq l_0\sum_{j \in V''} x_{ij} \hspace{3mm} \forall i \in D_0
\end{equation}

\begin{equation}
\label{eq.14}
\xi_i\leq Q \hspace{3mm} \forall i \in F
\end{equation}

\begin{equation}
\label{eq.15}
y_i\leq \xi_i \hspace{3mm} \forall i \in F
\end{equation}

\begin{equation}
\label{eq.16}
y_i\leq Q \hspace{3mm} \forall i \in D_0
\end{equation}

\begin{equation}
\label{eq.18}
x_{ij}\in \{0,1\} \hspace{3mm} \forall i \in V', \forall j \in V'', i \neq j
\end{equation}

The objective function (\ref{eq.objective_function}), to be minimized, consists of two components: the first one represents the total number of EVs used while, the second one measures the total time spent by the EVs outside the depot. Since our aim is firstly to minimize the EVs used and then the total route duration, the second component is normalized with regard to $l_0|N||D_0|$.

The group (\ref{eq.1}) of constraints imposes that each customer is 
visited exactly once while the group (\ref{eq.2}) guarantees that each clone of the RSs is used at most once. 
Moreover, constraints (\ref{eq.3DepIniziale}) and (\ref{eq.3}) assure that each clone of the depot (belonging to either $D_0$ or $D_{0'}$) is used at most once.
The flow conservation constraints, for each customer/RS, are imposed in  (\ref{eq.4}).

Constraints (\ref{eq.5_1}) assure that the arrival time to the vertex $j$ (i.e., either a customer or a final depot) is at least equal to the sum of the starting time from the previous vertex $i$ (i.e., either a customer or an initial depot), evaluated as $\tau_{i}+s_i$ and of the travel time to reach $j$ (i.e., $t_{ij}$). Similarly, constraints (\ref{eq.6_1}) guarantee that the arrival time to the vertex $j$ (i.e., either a customer or a final depot) is at least equal to the sum of the starting time from the previous RS $i$, expressed as $\tau_{i}+g(\xi_i-y_i)$ and of the travel time to reach $j$ (i.e., $t_{ij}$).
The group of constraints (\ref{eq.7}) assures that each customer and each initial/final depot is served inside its time window. Conditions (\ref{eq.9})-(\ref{eq.11}) assure that the vehicle capacity is never violated, considering that each customer request is indeed a request of pickup.

The group of constraints (\ref{eq.12}) assures that the vehicle capacity $C$ of each vehicle is never exceeded, considering that each customer request is of pickup and that $q_i=0 \hspace{3mm} \forall i \in F$.

The constraints (\ref{eq.13_bis})-(\ref{eq.13_ter}) logically link the arrival times to the routing variables. In particular, in (\ref{eq.13_bis}), if an arc $(i,j)$ is not traversed then the arrival time to the vertex $j$ has to be equal to zero. On the contrary, if the arc $(i,j)$ is traversed then the maximum arrival time at a vertex can be equal to the maximum route duration $l_0$. The same considerations can be applied to the group of constraints (\ref{eq.13_ter}).

Finally, the last groups of constraints (\ref{eq.14})-(\ref{eq.15}) assure that, for each EV, both the remaining battery level and the battery level reached at each RS are less or equal to the battery capacity $Q$. While, the constraints (\ref{eq.16}) impose that the remaining battery level of each EV at the final depot does never exceed $Q$.

It is worth noting that some preliminary evaluations have been carried out in a pre-processing phase at the aim of removing useless arcs from the graph $G$ modeling the problem. In particular, we remove the arcs $(i,j)$ violating the following condition:
\begin{equation}
\label{eq.19}
q_i + q_j \leq C \hspace{3mm} \forall i,j \in N, i \neq j
\end{equation}
since in this case the sum of the demands of customers $i$ and $j$ exceeds the vehicle capacity $C$.

\section{Solving the Time-Effective E-VRPTW by a Variable Neighborhood Search Branching}
\label{vns}
The proposed formulation of the E-VRPTW can be seen as a special case of $0-1$ Mixed Integer Linear Programming ($0-1$ MILP). Therefore, it is possible to adopt $0-1$ MILP solution methods, such as matheuristics.
Matheuristics are heuristic algorithms based on the solution of mathematical programming models.
Recently, their use for solving combinatorial optimization problems has allowed improving the state of the art on several problems such as also the ones coming from the real world (see \cite{Ball} for a survey on matheuristics).
Among them, we chose to apply the Variable Neighborhood Search Branching (VNSB) introduced in the seminal work of \cite{Hansen} and extended in \cite{Lazic}. It consists in adding linear constraints to the original problem for systematically changing the neighborhoods following the rules of the general Variable Neighborhood Search (VNS) schema. The idea of implementing a local search by adding to the MILP model a linear constraint, the so called {\it local branching constraint}, modeling the Hamming distance based neighborhood, has been firstly presented in the work of \cite{Fischetti} and is known as Local Branching Method. Therefore, the VNSB combines the VNS approach with the Local Branching one. To the best of our knowledge (\cite{Archetti,Doerner,Hansenb}), it has never been applied to the VRPTW before.

It is worth noting that when the VNSB is applied to a VRP not all the values of the right hand side ({\it rhs}) of the local branching constraint (i.e., the sizes of the Hamming distance based neighborhood) are feasible.
Indeed the minimum value of {\it rhs} that makes feasible a local branching constraint is $rhs=3$, since we may generate feasible solutions having Hamming distance three from a given feasible solution, in one of the following ways (see Figure \ref{cases}):

\begin{figure}[!h]   
\begin{center}
\begin{tabular}{c}
\includegraphics[width=15 cm]{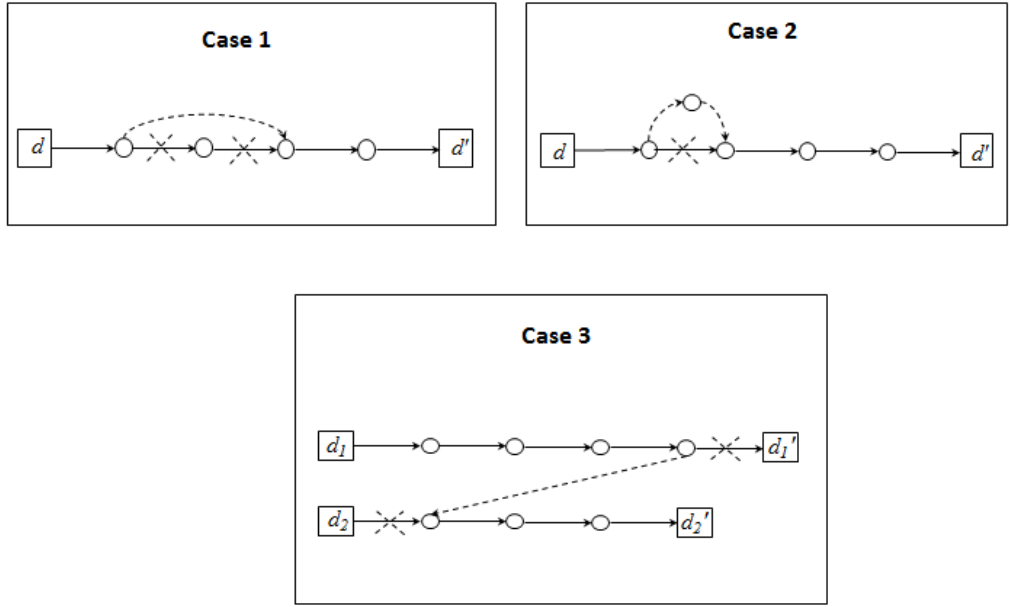}
\end{tabular}
\caption{Possible feasible solutions with Hamming distance three, where $d$, $d_1$ and $d_2$ are the clones of the depot used at the beginning of the route while $d'$, $d_1'$ and $d_2'$ are the ones used at the end} \label{cases}
\end{center}
\end{figure}

\begin{itemize}
\item \textbf{Case 1}: eliminating two consecutive arcs in a route and adding the arc linking the tail of the first one with the head of the second one. Note that since in this way we are eliminating a node from the route, the new route can be feasible only if such a node is an RS.
\item \textbf{Case 2}: eliminating one arc in a route and adding one arc linking the tail of the arc eliminated with a new node,
and another arc linking the new node with the head of the arc eliminated. Note that the new node can be only an RS, since the starting solution is feasible and then all the customer nodes are already covered.
\item \textbf{Case 3}: eliminating the last arc of a route, the first arc of another route and adding the arc linking the tail of the former with the head of the latter. Note that this corresponds to merge two routes.
\end{itemize}

On the other hand it is easy to see that any solution having Hamming distance less than three i.e., any transformation of a feasible solution into a solution that differs at most of two arcs is topologically infeasible.

With reference to the pseudo-code described in \cite{Hansen}, we set the parameter $k_{max}$, aimed to control the depth of the neighborhood explored during the shaking phase, equal to $10$ while the node time limit equal to five seconds.

\section{Computational Results}
\label{results}
In this section, we describe the computational results carried out on some benchmark instances taken from the literature. The MILP formulations have been implemented in AMPL and solved with the state of the art solver CPLEX 12.6 on a PC Intel Core i7, 3.20 GHz with 6GB RAM. The benchmark instances, available
at http://evrptw.wiwi.uni-frankfurt.de and described in \cite{Schneider}, have been generated from the ones proposed by Solomon for the VRPTW \cite{Solomon}.
\\In particular, three classes of instances, with $5$, $10$ and $15$ customers, respectively, have been experimented. In each of these classes, the instances have been distinguished according to the geographical distribution of the customers and the prefixes $R$, $C$ and $RC$ denote a random, clustered and mixed distribution, respectively. Moreover, each of these sub-classes has an additional classification where the prefixes $R1$, $C1$ and $RC1$ refer to a short scheduling horizon while $R2$, $C2$ and $RC2$, to a long scheduling horizon.
\\The aim of this section is twofold. Firstly, we compare the performances of the mathematical model proposed in \cite{Schneider} (hereafter, named as \textit{Model 1}) with ours (hereafter, named as \textit{Model 2}). This comparison is done considering:
\begin{itemize}
\item the number of vehicles used (denoted as $\mu_1$ and $\mu_2$, for \textit{Model 1} and \textit{Model 2}, respectively) that represents a common objective; \item the percentage gap of the Total Travel Distance (denoted as $\Delta TTD$); \item the percentage gap of the Total Time Outside the Depot (denoted as $\Delta TTOD$); \item the computational time required by both \textit{Model 1} and \textit{Model 2} (denoted as $CPU 1$ and $CPU 2$, respectively).
\end{itemize}
In particular, $\Delta TTD$ is computed as shown in the following:

\begin{equation}
\label{eq.TTD}
\Delta TTD=\frac{(TTD^{Model 2}-TTD^{Model 1})}{TTD^{Model 1}}\cdotp 100
\end{equation}

\noindent where $TTD^{Model 1}$ and $TTD^{Model 2}$ denote the total travel distance evaluated on the solution found by $Model 1$ and $Model 2$, respectively. Similarly, we define the $\Delta TTOD$ in the following way:

\begin{equation}
\label{eq.TTOD}
\Delta TTOD=\frac{(TTOD^{Model 2}-TTOD^{Model 1})}{TTOD^{Model 1}}\cdotp 100
\end{equation}

\noindent where $TTOD^{Model 1}$ and $TTOD^{Model 2}$ denote the total time spent outside the depot in the solution found by $Model 1$ and $Model 2$, respectively.\\
These numerical comparisons are shown in Table \ref{table.1} where the best results are highlighted in boldface. It is worth noting that all the computational results have been obtained with a CPU time limit set to $7,200.00$ seconds.
\begin{table}[htbp]
  \centering
\caption{Numerical Comparisons between \textit{Model 1} and \textit{Model 2}}
\label{table.1}
\scalebox{0.8}{
\begin{tabular}{|c|c|c|c|c|c|c|}
\hline
\textbf{Instance}&\textbf{$\mu_1$}&\textbf{$\mu_2$}&\textbf{$\Delta TTD$}&\textbf{$\Delta TTOD$}&\textbf{$CPU 1$}&\textbf{$CPU 2$}\\
\hline
C101C5&2&2&11.28&-12.83&0.53&3.06\\
C103C5&1&1&\textbf{0.00}&-12.87&0.27&0.16\\
C206C5&1&1&2.99&-16.44&3.70&94.94\\
C208C5&1&1&3.70&-60.89&0.89&4.15\\
R104C5&2&2&\textbf{0.00}&-2.18&1.00&20.68\\
R105C5&2&2&\textbf{0.00}&-7.95&0.36&2.09\\
R202C5&1&1&\textbf{0.00}&-13.78&2.15&19.30\\
R203C5&1&1&\textbf{0.00}&-53.36&5.64&736.64\\
RC105C5&2&2&\textbf{-3.11}&-25.47&8.80&58.31\\
RC108C5&2&2&\textbf{0.00}&-11.49&6.35&156.92\\
RC204C5&1&1&\textbf{0.00}&-66.06&84.71&\dag\\
RC208C5&1&1&\textbf{0.00}&-69.63&4.13&146.10\\
\hline
\hline
\textbf{Average}&&&\textbf{1.24}&\textbf{-30.92}&\textbf{10.73}&\textbf{767.25}\\
\hline
C101C10&3&3&21.37&-15.20&668.11&\dag\\
C104C10&2&2&18.03&-23.67&3625.46&\dag\\
C202C10&1&1&47.63&-2.18&3600.58&\dag\\
C205C10&2&2&28.43&-32.38&17.15&619.28\\
R102C10&3&3&\textbf{0.00}&-4.66&19.97&6,563.57\\
R103C10&2&2&0.10&-7.79&4721.45&\dag\\
R201C10&1&1&9.16&-12.72&122.76&\dag\\
R203C10&1&1&42.61&-45.93&3624.85&\dag\\
RC102C10&4&4&\textbf{0.00}&-11.01&49.18&\dag\\
RC108C10&3&3&5.03&-10.91&5214.34&\dag\\
RC201C10&1&1&14.91&-2.87&78.34&3,535.97\\
RC205C10&2&2&25.87&-8.96&90.96&\dag\\
\hline
\hline
\textbf{Average}&&&\textbf{17.76}&\textbf{-14.85}&\textbf{1,819.43}&\textbf{6,307.77}\\
\hline
C103C15&3&--&--&--&\dag&\dag\\
C106C15&3&3&39.06&-33.10&928.45&\dag\\
C202C15&2&2&52.43&-17.77&\dag&\dag\\
C208C15&2&2&56.37&-45.74&\dag&\dag\\
R102C15&5&--&--&--&\dag&\dag\\
R105C15&4&4&33.02&0.41&\dag&\dag\\
R202C15&2&2&63.96&-28.92&\dag&\dag\\
R209C15&1&1&39.62&-10.98&\dag&\dag\\
RC103C15&4&--&--&--&\dag&\dag\\
RC108C15&3&--&--&--&\dag&\dag\\
RC202C15&2&--&--&--&\dag&\dag\\
RC204C15&1&--&--&--&\dag&\dag\\
\hline
\hline
\textbf{Average}&&&\textbf{47.41}&\textbf{-22.68}&\textbf{6,154.74}&\textbf{7,200.00}\\
\hline
\end{tabular}
}
\end{table}

Concerning the number of vehicles used, in the sets with both $5$ and $10$ customers, \textit{Model 1} and \textit{Model 2} always give the same result. While, in the set with $15$ customers, this number is the same only in the cases in which \textit{Model 2} is able to detect a solution within the CPU time limit of $7,200.00$ seconds.

Moreover, the average percentage worsening of our model on $TTD$, considering the instance set with $5$ customers, is less than the one of \textit{Model 1} on $TTOD$ ($1.24\%$ against $30.92\%$, on average). It is worth noting that in $8$ of these instances (i.e., in the $66\%$ of this testbed), our model also optimizes $TTD$ (cases emphasized in boldface in Table \ref{table.1}). Considering, instead, the instance set with $10$ customers, the average percentage worsening of our model on $TTD$ is $17.76\%$ against $14.85\%$ on $TTOD$ provided by \textit{Model 1}. Moreover, in two cases, our model is suitable to also optimize $TTD$ (cases emphasized in boldface in Table \ref{table.1}). Finally, in the instance set with $15$ customers, for the cases in which \textit{Model 2} is suitable to find a solution within the CPU time limit, its average percentage worsening on $TTD$ is $47.41\%$ against $22.68\%$ of \textit{Model 1} on $TTOD$.

It is worth noting that for six of these instances with $15$ customers, our model is not suitable to find even a feasible solution within the CPU time limit, remarking the fact that the new version of the problem is more challenging than the original one.

Regarding the computational times, on average, our model requires $767.25$ seconds against $10.73$ seconds of \textit{Model 1}, on the instances with $5$ customers; $6,307.77$ seconds against $1,819.43$ seconds, on the instances with $10$ customers; $7,200.00$ seconds against $6,154.74$ seconds, on the instances with $15$ customers. This again remarks the more challenging nature of our problem.

It is finally worth remarking that in $11$ instances (i.e., in $30\%$ of the cases), our model determines solutions dominating those found by \textit{Model 1}, highlighted in boldface in Table \ref{table.1}.

However, in order to overcome the drawback due to the computational effort, we design a VNSB, as described in Section \ref{vns}. In particular, in Table \ref{table.2}, the results of our model are compared with the VNSB, setting the CPU time limit to $7,200.00$ seconds (as done for the MILPs).

In this case, the headers $\mu_3$ and $CPU_3$ denote, for each instance, the number of vehicles and the computational time required by the VNSB, respectively. Moreover, the headers $TTOD_2$ and $TTOD_3$ indicate, for each instance, the total time spent by the vehicles used outside the depot, provided by the proposed MILP formulation and the VNSB, respectively.

\begin{table}[htbp]
  \centering
\caption{Numerical Comparisons between \textit{Model 2} and the VNSB}
\label{table.2}
\scalebox{0.8}{
\begin{tabular}{|c|c|c|c|c|c|c|c|}
\hline
\textbf{Instance}&\textbf{$\mu_2$}&\textbf{$\mu_3$}&\textbf{$TTOD_2$}&\textbf{$TTOD_3$}&\textbf{$CPU_2$}&\textbf{$CPU_3$}\\
\hline
C101C5&2&2&1,262.84&1,262.84&3.06&0.70\\
C103C5&1&1&987.87&987.87&0.16&0.29\\
C206C5&1&1&1,296.82&1,296.82&94.94&14.16\\
C208C5&1&1&984.80&984.80&4.15&1.31\\
R104C5&2&2&196.17&196.17&20.68&2.43\\
R105C5&2&2&231.59&231.59&2.09&0.33\\
R202C5&1&1&234.16&234.16&19.30&1.62\\
R203C5&1&1&287.09&287.09&736.64&1.81\\
RC105C5&2&2&314.31&314.31&58.31&94.01\\
RC108C5&2&2&342.32&342.32&156.92&1.75\\
RC204C5&1&1&264.86&264.86&\dag&49.16\\
RC208C5&1&1&253.17&253.17&146.10&8.00\\
\hline
\hline
\textbf{Average}&&&\textbf{554.67}&\textbf{554.67}&\textbf{703.57}&\textbf{14.63}\\
\hline
C101C10&3&3&2,411.08&\textbf{2,335.20}&\dag&183.43\\
C104C10&2&2&1,805.67&\textbf{1,611.72}&\dag&73.43\\
C202C10&1&1&2,949.77&2,949.75&\dag&5.62\\
C205C10&2&2&2,525.77&2,525.77&619.28&0.79\\
R102C10&3&3&443.62&443.63&6,563.57&33.96\\
R103C10&2&2&351.90&\textbf{347.70}&\dag&93.49\\
R201C10&1&1&536.29&536.38&\dag&21.70\\
R203C10&1&1&540.29&\textbf{527.68}&\dag&4.62\\
RC102C10&4&4&571.26&571.25&\dag&59.01\\
RC108C10&3&3&515.11&\textbf{493.23}&\dag&143.71\\
RC201C10&1&1&793.52&793.52&3,535.97&40.46\\
RC205C10&2&2&618.44&\textbf{611.66}&\dag&37.51\\
\hline
\hline
\textbf{Average}&&&\textbf{1,171.89}&\textbf{1,145.62}&\textbf{6,307.77}&\textbf{58.14}\\
\hline
C103C15&--&3&--&\textbf{5,254.93}&\dag&240.46\\
C106C15&3&3&2,356.91&\textbf{2,173.09}&\dag&106.30\\
C202C15&2&2&3,936.32&\textbf{3,664.12}&\dag&139.83\\
C208C15&2&2&3,305.25&\textbf{2,819.47}&\dag&591.57\\
R102C15&--&5&--&\textbf{1,241.80}&\dag&559.52\\
R105C15&4&4&697.13&\textbf{567.82}&\dag&487.02\\
R202C15&2&2&999.72&\textbf{839.74}&\dag&349.29\\
R209C15&1&1&771.94&\textbf{621.79}&\dag&330.38\\
RC103C15&--&4&--&\textbf{1,074.90}&\dag&251.55\\
RC108C15&--&3&--&\textbf{974.88}&\dag&756.61\\
RC202C15&--&2&--&\textbf{3,222.49}&\dag&444.54\\
RC204C15&--&1&--&\textbf{2,742.46}&\dag&617.22\\
\hline
\hline
\textbf{Average}&&&\textbf{2,011.21}&\textbf{1,781.01}&\textbf{7,200.00}&\textbf{334.06}\\
\hline
\end{tabular}
}
\end{table}

It is worth noting that the VNSB always outperforms \textit{Model 2} for that concerns the total computational times. In fact, it requires on average: $14.63$ seconds against the $703.57$ seconds of \textit{Model 2}, on the instances with $5$ customers; $58.14$ seconds against $6,307.77$ seconds, on the instances with $10$ customers and finally, $334.06$ seconds against $7,200.00$, on the instances with $15$ customers.

Moreover, in some cases (highlighted in boldface in Table \ref{table.2}), the designed VNSB outperforms \textit{Model 2} also in terms of the solution quality. In fact, it determines a better feasible solution than the one detected by our model (within the same CPU time limit), with an average improvement of $45.05\%$ on the instances with $10$ customers and of $230.21\%$ on the ones with $15$ customers for which the proposed MILP formulation obtains a feasible solution in $7,200.00$ seconds.

\section{Conclusions and future works}
\label{conclusions}
In this work, the problem of routing and scheduling a fleet of Electric Vehicles (EVs) on a road network in order to serve a set of customers was addressed from both the modeling and the methodological point of view. In particular, the problem was firstly modeled, on a directed graph, as a VRPTW for EVs. Here we gave the original contribution of modeling the concept of partial recharge.
Then, a VNSB was designed for solving real case studies in reasonable amounts of time. To the best of our knowledge, such a matheuristic has never been proposed for solving a VRPTW before.

Computational results carried out on a set of benchmark instances showed two significant insights.

On one hand, our model improved the $TTOD$ of $23.02\%$, on average, compared to model proposed by \cite{Schneider}, worsening the $TTD$ only of $17.08\%$. Moreover, in the $30\%$ of the instances, our model obtained solutions dominating the ones found by \cite{Schneider}.

On the other hand, the designed VNSB was suitable to overcome the drawback due to the computational effort required by the proposed MILP, taking $95.92$ seconds, on average, against $4,248.46$ seconds of the latter.

In $13$ instances (i.e., $36\%$ of the cases) the VNSB obtained also solutions that are better than those found by the proposed MILP when it reached the CPU time limit, with an average relative improvement of $4.22\%$.

Finally, the designed VNSB found a solution, in reasonable amount of time, also in the six cases with $15$ customers for which the MILP was not suitable to determine even a feasible solution.

Future works concern the extension of both the proposed MILP and the VNSB with the aim of including the so called ``regenerative breaking'' that allows the EVs to recuperate a percentage of their battery consume, along the descents. Moreover, we will intend to consider also several different technologies for recharging the EVs at the RSs (e.g., to perform faster recharges) as proposed in \cite{Righini}. Finally, numerical experiments will be also carried out on large scale real world problems.


\end{document}